\documentclass{amsart}
\usepackage[mathscr]{eucal}
\usepackage{amsmath,amssymb,hyperref}
\usepackage{graphics}

\newtheorem{thm}{Theorem}[section]
\newtheorem{THM}{Theorem}
\newtheorem{PROP}{Proposition}

\newtheorem{prop}[thm]{Proposition}

\newtheorem{lemma}[thm]{Lemma}

\newtheorem{remark}[thm]{Remark}
\newtheorem{example}[thm]{Example}

\input xy
\xyoption{all}

\begin{document}

\title[Transversely Projective Foliations]
{Transversely Projective Foliations on Surfaces: Existence of Normal
Forms and Presciption of the Monodromy}
\author[F. Loray  and J.V. Pereira ]
{ Frank LORAY$^1$ and  Jorge Vit\'{o}rio PEREIRA$^{2}$}
\address{\newline $1$ IRMAR, Campus de Beaulieu, 35042 Rennes Cedex, France\hfill\break
$2$ IMPA, Estrada Dona Castorina, 110, Horto, Rio de Janeiro,
Brasil} \email{frank.loray@univ-rennes1.fr, jvp@impa.br}

 \maketitle

\begin{abstract}
We introduce a notion of normal form for transversely projective
structures of singular foliations on
complex manifolds. Our first main result says that this normal form
exists and  is unique when ambient space is two-dimensional. From
this result one obtains a natural way to produce invariants for
transversely projective foliations on surfaces. Our second main
result says that on projective surfaces one can construct singular
transversely projective foliations with prescribed monodromy.
\end{abstract}
\setcounter{tocdepth}{1}

\section{Introduction and Statement of Results}\label{S:intro}

\subsection{Singular Transversely Projective Foliations} Classically a
smooth holomorphic  transversely projective   foliation on a complex
manifold $M$ is a codimension one smooth holomorphic foliation
locally induced by holomorphic  submersions on $\mathbb P^1_{\mathbb
C}$ and with transitions functions in $\mathrm{PSL}(2,\mathbb C)$.
Among a  number of equivalent definitions that can be found in the
literature, we are particularly fund of the following one: $\mathcal
F$ is a {\bf transversely projective foliation} of a complex
manifold $M$ if there exists
\begin{enumerate}
\item  $\pi: P \to M$ a $\mathbb P^1$-bundle over $M$;
\item $\mathcal H$ a codimension one foliation of $P$ transversal to
the fibration $\pi$;
\item $\sigma: M \to P$ a holomorphic section transverse to
$\mathcal H$;
\end{enumerate}
such that $\mathcal F = \sigma^* \mathcal H$. The datum $\mathcal
P=(\pi:P \to M, \mathcal H, \sigma:M \to P)$ is the {\bf
transversely projective structure} of $\mathcal F$. A nice property
of this definition is that the isomorphism class of the $\mathbb
P^1$-bundle $P$ is an invariant canonically attached to the
foliation $\mathcal F$, whenever $\mathcal F$  has a leaf with
non-trivial holonomy, cf. \cite[page 177, Ex. 3.24.i]{Godbillon}.

In the holomorphic category the existence of smooth holomorphic
foliations imposes strong restrictions on the complex manifold. For
instance there exists a complete classification of smooth
holomorphic foliation on  compact complex surfaces, cf.
\cite{Brunella} and references there within. An interesting
corollary of this classification is that a rational surface carries
a holomorphic foliation if, and only if,
 it is a Hirzebruch surface and the foliation is a rational
fibration.

On the other hand the so called Riccati foliations on compact
complex surfaces $S$ , i.e., the foliations which are transversal to
a generic fiber of a rational fibration, are examples of foliations
which are transversely projective when restricted to the open set of
$S$ where the transversality of $\mathcal F$ with the rational
fibration holds.

The problem of defining a {\it good} notion of singular transversely
projective foliation on compact complex manifolds naturally emerges.
A first idea  would be to consider singular holomorphic foliations
which are transversely projective on Zariski open subsets. Although
natural, the experience shows that such concept is not very
manageable: it is too permissive. With an eye on applications one is
lead to impose some kind of regularity at infinity. A natural
regularity condition was proposed by Sc\'{a}rdua in \cite{scardua}.
Loosely speaking, it is imposed that the {\it transversely
projective structure} is induced by a global meromorphic triple of
$1$-forms. The naturality of such definition
 has been confirmed by the recent works of Casale on the
extension of Singer's Theorem \cite{Casale} and of Malgrange on
Non-Linear Differential Galois Theory \cite{Malgrange,
CasaleFourier}.

At this work we will adopt a variant of the above mentioned
 definition which maintains the geometric flavor of the
definition of a smooth transversely projective foliation given at
the beginning of the introduction. For us, $\mathcal F$ is a {\bf
singular transversely projective foliation} if there exists
\begin{enumerate}
\item  $\pi: P \to M$ a $\mathbb P^1$-bundle over $M$;
\item $\mathcal H$ a codimension one singular holomorphic foliation
of $P$ transverse to
the generic fiber of  $\pi$;
\item $\sigma: M \dashrightarrow P$ a meromorphic
section generically transverse to
$\mathcal H$;
\end{enumerate}
such that $\mathcal F = \sigma^* \mathcal H$. Like in the regular
case we will call the datum $\mathcal P=(\pi: P \to M, \mathcal H,
\sigma: M \dashrightarrow P)$ a {\bf singular transversely
projective structure} of $\mathcal F$.

A first remark is that unlike in the regular case the isomorphism
class of $P$  is not determined by $\mathcal F$ even when we suppose
that $\mathcal F$ is not {\it singular transversely affine}. In
general the $\mathbb P^1$-bundle $P$ is unique just up to
bimeromorphic bundle transformations. Thus the invariant that we
obtain is the bimeromorphism class of $P$. When $M$ is projective
this invariant is rather dull: any two $\mathbb P^1$-bundles over
$M$ are bimeromorphic.

To remedy this lack of unicity what we need is a
\subsection{Normal form for a singular transversely projective structure}
To a singular transversely projective structure $\mathcal
P=(\pi:P\to M,\mathcal H,\sigma)$ we associate the following objects
 on $M$:
\begin{itemize}
\item {\bf the branch locus},  denoted by $\text{Branch}(\mathcal P)$,
is the analytic subset of $M$ formed by the  points
$p\in M$ such that  $\sigma(p)$ is tangent to $\mathcal H$;
\item {\bf the indeterminacy locus}, denoted by  $\text{Ind}(\mathcal P)$,
is the analytic subset of $M$ corresponding to the  indeterminacy
locus of $\sigma$;
\item {\bf the polar divisor}, denoted by $(\mathcal
P)_{\infty}$, is the divisor on $M$ defined by the  direct image
under $\pi$ of the tangency divisor of $\mathcal H$ and the
one-dimensional foliation induced by the fibers of $\pi$.
\end{itemize}
Two transversely projective structures $\mathcal P=(\pi:P\to
M,\mathcal H,\sigma)$ and $\mathcal P'=(\pi:P'\to M,\mathcal
H',\sigma')$ are said to be {\bf bimeromorphically equivalent} if
there exists a bimeromorphism $\phi: P \dashrightarrow P'$ such that
$\phi^* \mathcal H'= \mathcal H$ and the diagram
\[
 \xymatrix{
 P \ar@{-->}[rr]^{\phi} \ar[dr]^{\pi} && P'
\ar[dl]_{\pi' } \\
& \ar@/^0.3cm/@{-->}[ul]^{\sigma} M \ar@/^-0.3cm/@{-->}[ur]_{\sigma'} } \\
\]
commutes.

We will say that a singular transversely projective structure
$\mathcal P$ is in {\bf normal form } when $\mathrm{cod} \, \,
\text{Branch}(\mathcal P) \le 2$ and the divisor $(\mathcal
P')_{\infty} - (\mathcal P)_{\infty}$ is effective, i.e.  $(\mathcal
P')_{\infty} - (\mathcal P)_{\infty}\ge0$, for every projective
structure $\mathcal P'$ bimeromorphic to $\mathcal P$ satisfying
$\mathrm{cod} \, \, \text{Branch}(\mathcal P') \le 2$.

\begin{THM}\label{T:normalform}
Let $\mathcal F$ be a singular transversely projective foliation on
a complex surface $S$. Every transversely projective structure
$\mathcal P$ of $\mathcal F$ is bimeromorphically equivalent to a
transversely projective structure in normal form. Moreover this
normal form is unique up to $\mathbb P^1$-bundle isomorphisms.
\end{THM}

We do not know if a normal form always exists  on higher dimensional
complex manifolds. Although when a normal forms  exists  our prove
of Theorem \ref{T:normalform} shows that it is unique.

From the unicity of the normal we can systematically produce
invariants for singular transversely projective foliations  on
complex surfaces. For singular transversely projective  foliations
on  the projective plane we define the

\subsection{Eccentricity of a Singular Transversely Projective Structure}

Let  $\mathcal P=(\pi:P \dashrightarrow \mathbb P^2, \mathcal H,
\sigma:\mathbb P^2 \dashrightarrow P)$ be a singular transversely
projective structure in normal form of  a foliation $\mathcal F$ of
the projective plane $\mathbb P^2$. We define the {\bf eccentricity}
of $\mathcal P$, denoted by $\mathrm{ecc}(\mathcal P)$, as follows:
if $L \subset \mathbb P^2$ is a generic line and  $P\vert_L$ is the
restriction of the $\mathbb P^1$-bundle $P$ to $L$ then we set
$\mathrm{ecc}(\mathcal P)$ as minus the self-intersection in
$P\vert_L$ of $\overline{\sigma(L)}$.

It turns out that the eccentricity of $\mathcal P$ can be easily
computed once we know the degree of the polar divisor. More
precisely we have the

\begin{PROP}\label{P:ecc}
Let $\mathcal F$ be a foliation on $\mathbb P^2$ and $\mathcal P$ a
singular transversely projective structure for $\mathcal F$ in
normal form. Then
\[
   \mathrm{ecc}(\mathcal P) =     \deg (\mathcal P)_{\infty} -  (\deg  (\mathcal F)
   +
    2) \, .
\]
\end{PROP}

We do not know if it is possible to give upper bounds for
$\mathrm{ecc}(\mathcal P)$ just in function of the degree of
$\mathcal F$. A positive result on this direction would be relevant
for what is nowadays called the Poincar\'{e} Problem.

The next result shows that $\mathrm{ecc}(\mathcal P)$ captures
dynamical information about $\mathcal F$ in some special cases.

\begin{PROP}\label{P:dynamics}
Let $\mathcal F$ be a quasi-minimal singular transversely projective
foliation of $\mathbb P^2$ and $\mathcal P$ be a transversely
projective structure for $\mathcal F$ in normal form. If the
monodromy representation of $\mathcal P$ is not minimal then
\[
\mathrm{ecc}(\mathcal P) > 0 \, .
\]
\end{PROP}

An immediate corollary  is that transversely projective structures
in normal form of Hilbert modular foliations on $\mathbb P^2$ have
positive eccentricity. This follows from Proposition
\ref{P:dynamics} and the well-known facts that these foliations are
transversely projective, quasi-minimal and with monodromy contained
in $\mathrm{PSL}(2,\mathbb R)$, cf. \cite[Theorem 1]{LG}.

\subsection{The Monodromy Representation}
A very important invariant of a projective structure $\mathcal
P=(\pi:P\to M, \mathcal H, \sigma:M\dashrightarrow P)$, is the {\bf
monodromy representation}. It is the representation of $\pi_1(M
\setminus \vert (\mathcal P)_{\infty} \vert)$ into
$\mathrm{PSL}(2,\mathbb C)$ obtained by lifting paths on $M
\setminus \vert (\mathcal P)_{\infty} \vert$ to the leaves of
$\mathcal H$.

Given a hypersurface $H \subset M$ and a representation $\rho:
\pi_1(M\setminus H )\to \mathrm{PSL}(2,\mathbb C)$, one might ask if
there exists a foliation $\mathcal F$ of $M$ with transversely
projective structure $\mathcal P$ whose monodromy is $\rho$.

We will show in \S \ref{S:obstrucao} that the answer is in general
no: there are local obstructions to solve the {\it realization
problem}.

On the other hand if the ambient is two-dimensional and the
representation $\rho$ lifts to a representation  $\tilde \rho:
\pi_1(M\setminus H) \to \mathrm{SL}(2,\mathbb C)$ then we have the

\begin{THM}\label{T:RH}
Let $S$ be a projective surface  and $H$ a reduced hypersurface on
$S$. If
$$\rho:\pi_1(S\setminus H ) \to \mathrm{PSL}(2,\mathbb C)$$ is a homomorphism
which lifts to a homomorphism $\tilde{\rho}:\pi_1(S\setminus  H )
\to \mathrm{SL}(2,\mathbb C)$ then there exists a singular
transversely projective foliation $\mathcal F$ with a   singular
transversely projective structure in normal form $\mathcal P$ such
that
\begin{enumerate}
\item $H- (\mathcal P)_{\infty}\ge 0 $;
\item $\rho$  is the monodromy representation
of $\mathcal P$;
\item If $\rho$ is not solvable then $\mathcal F$ admits a unique singular transversely projective structure in normal form.
\end{enumerate}
\end{THM}

We point out that the result (and the  proof here presented)  holds
for higher dimensional projective manifolds if one supposes that $H$
is a normal crossing divisor, cf. \S \ref{S:RH} for details.





\section{Generalities}\label{S:singular}

\subsection{A local description of $\mathcal H$}
Let $\Delta^n \subset \mathbb C^n$ be a polydisc and $\pi:P\to
\Delta^n$ be a $\mathbb P^1$-bundle. Since the polydisc is a Stein
contractible space we can suppose that $P$ is the projectivization
of the trivial rank $2$ vector bundle over $\Delta^n$ and write
$\pi(x,[z_1:z_2])=x$. If $\mathcal H$ is a codimension one foliation
of $P$ generically transversal to the fibers of $\pi$ then
$\pi^*\mathcal H$ is induced by a $1$-form $\Omega$ that can be
written as
\[
  \Omega = z_1dz_2 - z_2 dz_1 + \alpha z_1^2 + \beta z_1\cdot
  z_2 + \gamma z_2^2 \, ,
\]
where   $\alpha, \beta$ and $\gamma$ are meromorphic $1$-forms on
$\Delta^n$.  The integrability  condition $\Omega \wedge d\Omega =
0$ translates into the relations
\begin{equation}\label{E:triplet}
\left\{\begin{matrix}
d\alpha &=& \hfill\alpha \wedge \beta \\
d\beta &=& 2\alpha \wedge \gamma \\
d\gamma &=& \hfill\beta \wedge \gamma
\end{matrix} \right.
\end{equation}

 The divisor of poles of $\Omega$ corresponds to the fibers of
$\pi$ that are tangent to $\mathcal H$, i.e., if $\mathcal C$
denotes the $1$-dimension foliation induced by the fibration $\pi$
then
\[
   (\Omega)_{\infty} = \mathrm{tang}(\mathcal H, \mathcal C) \, .
\]

Associated to $\Omega$ we have an  integrable differential
$\mathfrak{sl}(2,\mathbb C)$-system on the trivial rank $2$ vector bundle over
$\Delta^n$ defined by
\[
dZ=A\cdot Z\ \ \ \text{where}\ \
\
A=\left(\begin{matrix} -{\frac{\beta}{ 2}}&-\gamma\\
\alpha&{\frac{\beta}{ 2}}
\end{matrix}\right)\ \ \ \text{and}\ \ \
Z=\left(\begin{matrix} z_1\\z_2\end{matrix}\right)
\]
The matrix $A$ can  be thought as a meromorphic differential
$1$-form on $\Delta^n$ taking values in the Lie algebra
$\mathfrak{sl}(2,\mathbb C)$ and satisfying the integrability
condition $dA+A\wedge A=0$. Darboux's Theorem (see \cite{Godbillon},
III, 2.8, iv, p.230) asserts that  on any simply connected open
subset $U\subset \Delta^n \setminus (\Omega)_{\infty}$ there exists
a holomorphic map
$$\Phi:U\to \mathrm{SL}(2,\mathbb C)\ \ \ \text{such that}\ \ \ A=\Phi^* M$$
where $ M$ is the Maurer-Cartan $1$-form on $\mathrm{SL}(2,\mathbb
C)$. Moreover, the map $\Phi$ is unique up to a left  composition
with an element in  $\mathrm{SL}(2,\mathbb C)$. For every  $v \in
\mathbb C^2$  the sections
\begin{eqnarray*}
   \varphi_v: U &\to& U \times \mathbb C^2 \\
              x &\mapsto& (x,\Phi(x)\cdot v)
\end{eqnarray*}
are  solutions of the differential system above. It follows that the
application
\begin{eqnarray*}
\phi: U \times \mathbb P^1 &\to& U \times \mathbb P^1 \\
(x,[z_1,z_2]) & \mapsto & ( x, [\Phi(x)(z_1,z_2)] ) \, ,
\end{eqnarray*}
conjugates the foliation $\mathcal H\vert_U$ with the one induced by
the submersion $U\times \mathbb P^1 \to \mathbb P^1$.

We have just described  $\mathcal H$ over the points outside
$(\Omega)_{\infty}$. Now we turn our attention to

\subsection{The behaviour of  $\mathcal H$ over a generic point of
$(\Omega)_{\infty}$}\label{S:generic}

Let $W$ be an analytic subset of the support of $(\Omega)_{\infty}$.
We will set $S(W)$ as
\[
  S(W)=\pi^{-1}(W) \cap \mathrm{sing}(\mathcal H) \, .
\]

We will  start by analyzing $\mathcal H$ over the irreducible
components $H$ of $(\Omega)_{\infty}$ for which $\pi^{-1}(H)$ is
$\mathcal H$-invariant.

\begin{lemma}\label{L:trivializa}
Let  $H$ be an irreducible component of  the support of
$(\Omega)_{\infty}$ and
\[
   V = \{ p \in H \, \, \text{such that}  \, \, \pi^{-1}(p) \subsetneq \mathrm{sing}(\mathcal
   H) \text{ and } H \text{ is smooth at } p  \}.
\]
If $\pi^{-1}(H)$ is $\mathcal H$-invariant then for every $p \in V$
there exists a neighborhood $U$ of $p$, a map $\varphi:(U,H\cap U)
\to (\mathbb C,0)$ and
 a Riccati foliation $\mathcal R$ on $(\mathbb C,0) \times
\mathbb P^1$ such that $$\mathcal H = \varphi^* \mathcal R.$$ In
particular $\pi_{|S(V)} : S(V) \to V$ is an \'{e}tale covering of $V$ of
degree $1$ or $2$.
\end{lemma}
\begin{proof}
Since $\mathrm{cod } \, \mathrm{sing}(\mathcal H)\le 2$ then $V$ is
a dense open subset of $H$.

Let $p \in V$ and $F\in \mathcal O_{\Delta^n,p}$ be a local equation
around for the poles of $\Omega$. Since $p \in V$ at least one of
the holomorphic $1$-forms $F\alpha,F\beta,F\gamma$ is non-zero at
$p$. After applying a  change of coordinates of the form
\[
(x,[z_1:z_2]) \mapsto (x,[a_{11}  z_1 + a_{12}z_2 : a_{21}  z_1 +
a_{22}z_2])
\]
where \[  \left(
                               \begin{array}{cc}
                                 a_{11} & a_{12} \\
                                 a_{21} & a_{22} \\
                               \end{array}
                             \right) \in \mathrm{GL}(2,\mathbb C), \]
                              we can assume that $F\alpha,F\beta$  and $F\gamma$ are
non-zero at $p$.

From the relation $d\alpha =  \alpha \wedge \beta $ we promptly see
that the holomorphic $1$-form $F\alpha$ besides being non-singular
is also integrable. It follows from Frobenius integrability Theorem
and the $\mathcal H$-invariance of $H$ that there exist a local
system of coordinates $(x,y_2,\ldots,y_n):U\to\mathbb C^n$ where $p$
is the origin of $\mathbb C^n$, $F=x^n$ for a suitable $n \in
\mathbb N$ and $F\alpha = h_0 dx$ for some $h_0 \in \mathcal
O^*_{\Delta^n,p}$.

Again from the relation $d\alpha  = \alpha  \wedge \beta$ and the
fact that $F\beta(p) =(x^n \beta)(p) \neq 0$ it follows that there
exists $h_1 \in \mathcal O^*_{\Delta^n,p}$ such that
\[
  \beta =   -\frac{dh_0}{h_0} + h_1 \frac{dx}{x^n} .
\]

After performing the holomorphic change of variables
$$
       (x,[z_1:z_2]) \mapsto \left(x, [ {h_0  z_1} :{z_2 + (1/2)h_1  z_1    }]\right)
$$
we can suppose that $(\alpha, \beta)=(\frac{dx}{x^n}, 0)$.

The conditions $d\beta=2\alpha \wedge \gamma$ and $d\gamma = \beta
\wedge \gamma$ imply that $\gamma$ depends only on $x$: $\gamma=
b(x)\frac{dx}{x^n}$, with $b$ holomorphic. Note that on this new
coordinate system we can no longer suppose that
$F\gamma(p)=x^n\gamma(p) \neq 0$.

Thus on this  new coordinate system
\[
\Omega = z_1dz_2 - z_2 dz_1 + z_1^2 \frac{dx}{x^n}  + z_2^2 b(x)
\frac{dx}{x^n} \, .
\]
It follows that on $\pi^{-1}(q)$, $q \in V$, we  have one or two
singularities of $\mathcal H$: one when $b(0)=0$ and two otherwise.
\end{proof}

\noindent {\bf A word about the terminology:} Further on when we
refer to the {\bf transverse type} of an irreducible curve of
singularities we will be making reference to the type of singularity
of the associated Riccati equation given by the above proposition.

\medskip

Let us now analyze $\mathcal H$ over the irreducible components $H$
of $(\Omega)_{\infty}$ for which $\pi^{-1}(H)$ is not $\mathcal
H$-invariant. In  the notation of lemma \ref{L:trivializa} we have
the

\begin{lemma}\label{L:correcao}
If $\pi^{-1}(H)$ is not $\mathcal H$-invariant then  $\pi_{|S(V)} :
S(V) \to V$ admits an unique holomorphic section.
\end{lemma}
\begin{proof} Let $p \in V$ be an arbitrary point. Without loss of generality we
can assume that $F\alpha,F\beta$ and $F\gamma$ are non-zero at $p$
and that $H$ is not invariant by the foliation induced by $\alpha$,
cf. proof of lemma \ref{L:trivializa}.

Assume also that  $\ker \alpha(p)$ is transverse to $H$. Thus there
exists a suitable local coordinate system $(x,y,y_3,\ldots, y_n):U
\to \mathbb C^n$ where  $p$ is the origin, $F=x^n$ for some  $n \in
\mathbb N$ and $F \alpha = h_0 dy$ for some $h_0 \in \mathcal
O^*_{\Delta^n,p}$.

The condition $d\alpha=\alpha\wedge\beta$ implies that $\beta=n
{\frac{dx}{ x}}+h_1\cdot \alpha$ with $h_1$ meromorphic at $p$.
Since $F\beta=x^n\beta$ is holomorphic and does not vanish at $p$
the same holds for $h_1$, i.e., $h_1 \in \mathcal O^*_{\Delta^n,p}$.
Thus if we apply the holomorphic change of coordinates
\[
 ((x,y,y_3,\ldots, y_n), [z_1:z_2] \mapsto ((x,y,y_3,\ldots,
 y_n),[z_1: h_0\cdot z_2 + h_1 \cdot z_1])
\]
we have $d\beta=0$.

Combining $0=d\beta = 2\alpha \wedge \gamma$ with $d\gamma=\beta
\wedge \gamma$ we deduce that $\gamma =x^n h_3(y)\alpha$ for some
meromorphic function $h_3$. Since $\gamma$ has poles contained in
$H=\{x=0\}$, $h_3$ is in fact holomorphic and consequently $\mathcal
H$ is induced by the $1$-form
\begin{equation}\label{E:fake}
x^n(z_1 dz_2 - z_2 dz_1) + (dy) z_1^2 + (x^n h_3(y) dy)z_2^2 \, .
\end{equation}
It is now clear that the singular set of $\mathcal H$ is given by
$\{x=0\} \cap \{z_1=0\}$. Thus there exists an open subset
$V_0\subset V$  for which $S(V_0)$ is isomorphic to $V_0$. Since
$S(V)$ does not contain fibers of $\pi_{|S(V)}$ this is sufficient
to prove the lemma.
\end{proof}

\begin{remark}\rm\label{R:util}
These irreducible components of $(\Omega)_{\infty}$ are a kind of
{\it fake} or {\it apparent} singular set for the transversely
projective structures. More precisely, after the fibred  birational
change of coordinates
\[
  ((x,y,y_3,\ldots, y_n), [z_1:z_2] \mapsto ((x,y,y_3,\ldots,
 y_n),[z_1: x^n z_2])
\]
the foliation  induced by (\ref{E:fake}) is completely transversal
to the fibres of the $\mathbb P^1$-fibration and has a product
structure as in the case $H$ is $\mathcal H$-invariant.
\end{remark}

\subsection{Elementary Transformations}\label{S:elementary}

Still in the local setup, let $H$ be a smooth and  irreducible
component of the support of $(\Omega)_\infty$ and let $S \subset
\pi^{-1}(H)$ be a holomorphic section of the restriction of the
$\mathbb P^1$-bundle over $M$. An {\bf elementary transformation}
$\mathrm{elm}_S: \Delta^n \times \mathbb P^1 \dashrightarrow
\Delta^n \times \mathbb P^1$ with center in $S$ can be described as
follows: first we blow-up $S$ on $M$ and then we contract the strict
transform of $\pi^{-1}(H)$. If $F=0$ is a reduced equation of $H$
and $S$ is the intersection of $H\times \mathbb P^1$ with the
hypersurface $z_2=0$ then $\mathrm{elm}_S$ can be explicitly written
as
\begin{eqnarray*}
\mathrm{elm}_S: \Delta^n \times \mathbb P^1 &\dashrightarrow&
\Delta^n \times
\mathbb P^1 \\
 ( x , [z_1:z_2] ) &\mapsto& (x, [  F(x) z_1  :  z_2 ] )
\end{eqnarray*}
modulo $\mathbb P^1$-bundle isomorphisms on the source and the
target.

We are interested in describing the foliation
$({\mathrm{elm}_S})_*\mathcal H=(\mathrm{elm}_S^{-1})^* \mathcal H$.
More specifically we want to understand how the divisors
$\mathrm{tang}(\mathcal H,\mathcal C)$ and
$\mathrm{tang}(e_*\mathcal H,\mathcal C)$ are related, where
$\mathcal C$ denotes the one dimension foliation induced by the
fibers to $\Delta^n \times \mathbb P^1 \to \Delta^n$. We point out
that the analysis we will now carry on can be found in the case
$n=1$ in \cite[pages 53--56]{Brunella}. The arguments that we will
use are essentially the same. We  decided to include them here
thinking on  readers' convenience.

Let $k$ be the order of $(\Omega)_{\infty}$ along $H$.  Since
$\mathrm{elm}_S^{-1}(x,[z_1:z_2])= (x,[z_1:F(x)z_2])$ it follows
that
\[
   (\mathrm{elm}_S^{-1})^* \Omega = F ( z_1dz_2 - z_2 dz_1) +  \alpha z_1^2 + F\left(\beta + \frac{dF}{F}\right) z_1\cdot
  z_2 + F^2\gamma z_2^2 \, .
\]
Thus the foliation $(\mathrm{elm}_S)_*\mathcal H$ is induced by the
meromorphic $1$-form
\[
\widetilde\Omega =( z_1dz_2 - z_2 dz_1) +  \frac{\alpha}{F} z_1^2 +
\left(\beta + \frac{dF}{F}\right) z_1\cdot
  z_2 + F\gamma z_2^2 \,
\]

In order to describe $(\Omega)_{\infty}$ we will consider three
mutually  exclusive cases:
\begin{enumerate}
\item \underline{$S$ is not contained in $\pi^{-1}(H) \cap
\mathrm{sing}(\mathcal H)$:} This equivalent to say that $F^k
\alpha$ is not identically zero when restricted to $H$. Therefore
\[
  (\widetilde \Omega)_{\infty} = (\Omega)_{\infty} + H \, .
\]
\item \underline{${S \subset \pi^{-1}(H) \cap
\mathrm{sing}(\mathcal H)}$ but $S$ is not equal to $\pi^{-1}(H)
\cap \mathrm{sing}(\mathcal H)$}: This corresponds to  $(F^k
\alpha)_{|H} \equiv 0$ while  $(F^k \beta)_{|H} \not \equiv 0$. If
$k\ge 2$ then it follows that
\[
  (\widetilde \Omega)_{\infty} = (\Omega)_{\infty} \, .
\]
When $k=1$ we have two possible behaviors
\[
  (\widetilde \Omega)_{\infty} = \left\lbrace \begin{array}{ll}
                                                (\Omega)_{\infty} -H & \text{ when }   \beta +  \frac{dF}{F} \text{ is holomorphic.} \\
                                                (\Omega)_{\infty}  &  \text{  otherwise .}
                                              \end{array} \right.
\]
\item \underline{$S$ is equal to
$\pi^{-1}(H) \cap \mathrm{sing}(\mathcal H)$}: Here  $(F^k
\alpha)\vert_{H} \equiv (F^k \beta)\vert_{H}  \equiv 0$ while $(F^k
\gamma)_{|H} \not \equiv 0$.When $k=1$ it follows that
\[
  (\widetilde \Omega)_{\infty} = (\Omega)_{\infty} \, .
\]
When $k\ge 2,$ if we set  $k'$ as  the smallest positive integer for
which
 $(F^{k'}\alpha)\vert_H \equiv (F^{k'+1}\beta)\vert_H\equiv
0 $ then
\[
  (\widetilde \Omega)_{\infty} = (\Omega)_{\infty} - (k-k') H .
\]
\end{enumerate}

In the global setup the picture is essentially the same, i.e., if
$\pi:P \to M$ is  a $\mathbb P^1$-bundle over $M$, $H$ is a smooth
hypersurface on $M$ and $s:M\to P\vert_H$ is a holomorphic section
then we build up a new $\mathbb P^1$-bundle by  blowing up the image
of $s$ in $P$ and  contracting the strict transform $\pi^{-1}(H)$.
The local analysis just made can be applied, as it is,  on the
global setup.

In the case $P$ is the projectivization  of a rank $2$ holomorphic
vector bundle over a complex manifold $M$ then the elementary
transformations just described are  projectivizations of the so
called elementary modifications, see \cite[pages 41--42]{Friedman}.

\section{Existence and Unicity of the Normal Form}

\subsection{Existence of a Normal Form I: A Particular  Case}Let $\mathcal
P=(\pi:P\to S,\mathcal H,\sigma:S \dashrightarrow P)$ be a
transversely projective structure for a foliation $\mathcal F$ on  a
complex surface $S$. We will now prove the existence of a normal
form for $\mathcal P$ under the additional assumptions that the
irreducible components of the support of $(\Omega)_{\infty}$ and the
 codimension one irreducible components of
$\mathrm{Branch}(\mathcal P)$ are smooth.

Let $H$ be an irreducible component of $(\mathcal P)_{\infty}$ of
multiplicity $k(H)$ and, as in \S \ref{S:generic}, let $S(H)$ be
given by $ S(H)=\pi^{-1}(H) \cap \mathrm{sing}(\mathcal H) \, .$
Thus ( see lemmata \ref{L:trivializa} and \ref{L:correcao} ) $S(H)$
is an analytic subset of $\pi^{-1}(H)$ formed by a finite union of
fibers together with a one or two-valued holomorphic section $s$ of
$P\vert_H$. Note that to assure that $s$ is in fact holomorphic, and
not just meromorphic, we have used that $H$  is a curve, i.e., we
have used that $S$ is a surface.

If $s$ is two-valued and $k(H)>1$ then the elementary transformation
centered in any of the branches of $s$ (we are, of course, restring
to a simply-connected open set where $s$ does not ramifies) will not
change the order of poles of $\Omega)$, i.e., the order of poles is
already minimal. This follows from the fact that the transverse type
of  $s$  is  reduced, it is in fact (cf. \cite[page 54]{Brunella}) a
saddle-node.

If $s$ is two-valued and $k(H)=1$ then $s$ does not ramifies. In
fact,  if the quotient of eigenvalues along   one of the branches of
$s$ is $\lambda$ then, by Camacho-Sad index theorem, the other
branch will have quotient of eigenvalues equal to $-\lambda$. Thus
ramification of $s$ leads to absurdity $\lambda \neq 0$ and
$\lambda=-\lambda$. If the quotient of eigenvalues  of the branches
of $s$ are  not integers then we are in a minimal situation. On the
contrary if the quotient of eigenvalues of one of the branches of
$s$, say $s_+$, is a positive integer, say $\lambda_+$, then by an
elementary transformation centered at $s_+$ we will obtain two new
sections of singularities one of them with transverse type
$\lambda_+ - 1$. After $\lambda_+$ successive elementary
transformations we will arrive at a transversely projective
structure, still denoted by $\mathcal P$, where $H$ does not belong
to the support of $(\mathcal P)_{\infty}$ (linearizable transverse
type) or $s$ is one-valued(Poincar\'{e}-Dulac transverse type).

If $s$ is one-valued, $k(H)=1$  and $H$ is $\mathcal H$-invariant
then an elementary transformation centered in $s$ will either
transform $\mathcal H$ to a foliation with $k(H)=1$ but now with $s$
two-valued. It changes the transverse type from saddle-node (with
weak separatrix in the direction of the fibration) to
Poincar\'{e}-Dulac. The important fact is that it does not changes
$k(H)$.

If $s$ is one-valued,  $k(H)>1$ and $H$ is $\mathcal H$-invariant
then we have two possibilities. The first is when the transverse
type is degenerated. An elementary transformation centered in $s$
will drop the multiplicity of $H$ on $(\mathcal P)_{\infty}$. The
second possibility is when the transverse type is nilpotent. On this
last case the multiplicity is stable by elementary transformations,
cf. \cite[page 55-56]{Brunella}.

If $s$ is one-valued and $H$ is not $\mathcal H$-invariant then an
elementary transformation centered in $s$ will decrease $k(H)$ by
one, compare with remark \ref{R:util}. Of course if $k(H)$ reaches
zero then the resulting foliation is smooth over a generic point of
$H$.

In resume after applying a finite number of  elementary
transformations we arrive at a projective structure, still denoted
by $\mathcal P$, for which $(\mathcal P)_{\infty}$ has  minimal
multiplicity in the same  bimeromorphic equivalence class. Note also
that no  codimension one components have been added to $\vert
(\mathcal P)_{\infty}\vert \cup \mathrm{Branch}(\mathcal P)$ along
the process.

Of course there are distinct  biholomorphic equivalence class of
projective structure with the same property. To rigidify we have to
consider $\mathrm{Branch}(\mathcal P)$.

Let now $H$ be an irreducible codimension one component of
$\mathrm{Branch}(\mathcal P)$. First suppose that $H$ is contained
in the support of $(\mathcal P)_{\infty}$. The restriction of
$\sigma$ to $\pi^{-1}(H)$ determines $s$ a natural candidate for
center of an elementary transformation. As before, keep the same
notation from the projective structure obtained after applying the
elementary transformation centered in $s$. Two things can happen:
(1) ${\sigma}\vert_{H} \subset \mathrm{sing}(\mathcal H)$; or (2)
${\sigma}\vert_{H} \not \subset \mathrm{sing}(\mathcal H)$. In case
(2) we are done. In case (1) we are in a situation no different from
the one that we started with. If we iterate the process and keep
falling in case (1) we deduce that $\sigma$ follows the infinitely
near singularities of $\mathcal H\vert_H$  and therefore must be an
$\mathcal H$-invariant hypersurface. Of course this is not the case
since in  the definition of a transversely projective structure we
demand that $\sigma$ is generically transverse to $\mathcal H$.

It remains to consider the case $H$ is not contained in the support
of $(\mathcal P)_{\infty}$. The elementary transformation centered
on $s$, the restriction of $\sigma$ to $\pi^{-1}(H)$, yields a
projective structure for which we have added $H$ with multiplicity
one in $(\mathcal P)_{\infty}$. So we have reduced to the case just
analyzed: $H$ is contained in the support of $(\mathcal
P)_{\infty}$.

In resume we have proved the

\begin{prop}\label{P:particular}
Let $\mathcal P=(\pi:P\to S,\mathcal H,\sigma:S \dashrightarrow P)$
be a transversely projective structure for a foliation $\mathcal F$
on a complex surface $S$. Suppose that  the irreducible components
of the support of $(\Omega)_{\infty}$ and the  codimension one
irreducible components of $\mathrm{Branch}(\mathcal P)$ are smooth.
Then there exists $\mathcal P'$ a transversely projective structure
in normal form bimeromorphically equivalent to  $\mathcal P$.
\end{prop}

Before dealing with the unicity of the normal form we will prove the

\subsection{Existence of a Normal Form II: The General Case}
To prove the existence of a normal form for a general transversely
projective structure $\mathcal P=(\pi:P\to S,\mathcal H,\sigma:S
\dashrightarrow P)$ for a foliation $\mathcal F$  we  proceed as
follows.

We start by taking an embedded resolution of the support of
$(\mathcal P)_{\infty}$ and of the codimension one components of
$\mathrm{Branch}(\mathcal P)$, i.e., we take a bimeromorphic
morphism $r:\tilde S \to S$ such that $r^*|(\mathcal P)_{\infty}|$
is a divisor with smooth irreducible components  and the codimension
one components of $r^*(\mathrm{Branch}(\mathcal P))$ are also
smooth. We will now work with $\widetilde{ \mathcal P}= r^*\mathcal
P$ a transversely projective structure for $\widetilde{ \mathcal
F}=r^*\mathcal F$.

Proposition \ref{P:particular} implies that there exists a
transversely projective structure $\widetilde{\mathcal P}'$ in
normal form bimeromorphic to $\widetilde{\mathcal P}$. If $U=\tilde
S \setminus D$, where $D$ denotes the exceptional divisor of $r$,
then $\mathcal P'= (r\vert_U)_* (\widetilde{\mathcal P}'\vert_U)$ is
a transversely projective structure in normal form for $\mathcal F$
defined on the complement of a finite number of points. It follows
from Hartog's extension Theorem that to extend $\mathcal P'$ it is
sufficient to extend the $\mathbb P^1$-bundle P'. To conclude we
have just to apply the following

\begin{lemma}\label{L:extensao}
Let $\widetilde \pi:\widetilde P \to \widetilde S$ be a $\mathbb
P^1$-bundle over a compact complex surface $\widetilde S$ and let
$r:\widetilde S \to S$ be a  bimeromorphic morphism with exceptional
divisor $D$. Then there exists a $\mathbb P^1$-bundle $\pi:P \to  S$
and a map $\phi:\widetilde P \to P$ such that $\phi\vert_{\widetilde
\pi^{-1}(\widetilde {S} \setminus D)}$ is a $\mathbb P^1$-bundle
isomorphism.
\end{lemma}
\begin{proof} Let $\widetilde U$ be a sufficiently small neighborhood of the
support of $D$. We can assume that $U=r(\widetilde U)$ is a Stein
subset of $S$. Outside $\widetilde U$ there is no problem at all:
the map $r$ is an biholomorphism when restricted to $\widetilde S
\setminus \widetilde U$.

Suppose first that there exists a rank 2 vector bundle $\tilde E$
over $\tilde U$ such that $\widetilde{P}\vert_{\tilde U} = \mathbb
P(\tilde E)$. If $\widetilde{\mathcal E}$ denotes the  sheaf of
sections of $\tilde E$, $U=r(\widetilde U)$ and $p=r(E)$ then
Grauert's direct image Theorem assures that $\phi_*
\widetilde{\mathcal E}$ is a coherent $\mathcal O_V$-sheaf. Moreover
$\phi_* \widetilde{\mathcal E}$ is locally free when restricted to
$U\setminus\{p\}$. If $\mathcal E^{\vee \vee} =
\mathrm{Hom}(\mathrm{Hom}(r_*\mathcal E,\mathcal O_V), \mathcal
O_V)$  then $\mathcal E^{\vee \vee}$ is a reflexive sheaf. Since we
are in dimension two $\mathcal E^{\vee \vee}$ is in fact locally
free, cf. \cite[Proposition 25, page 45]{Friedman}. Thus $\mathcal
E^{\vee \vee}$ is the sheaf of sections of some rank two vector
bundle $E$. This is sufficient to prove the lemma under the
assumption that $\widetilde{P}\vert_{\tilde U} = \mathbb P(\tilde
E)$. Now we will show that this is always the case.

The obstruction to a $\mathbb P^1$-bundle over $\tilde U$ be the
projectivization of a rank two vector bundles lies in $\mathrm
H^2(\tilde U,\mathcal O_{\tilde U}^*)$, cf. \cite[page 190]{BPV}. It
follows from the exponential sequence that  $\mathrm H^3(\tilde
U,\mathbb Z)= \mathrm H^2(\tilde U, \mathcal O_{\tilde U})=0$
implies that $\mathrm H^2(\tilde U,\mathcal O_{\tilde U}^*)=0$. But
$\tilde U$ has the same type of homotopy of a tree of rational
curves. Thus $\mathrm H^3(\tilde U,\mathbb Z)=0$. On the other hand,
by \cite[Theorem 9.1.(iii), pages 91--92]{BPV}, $\mathrm H^2(\tilde
U, \mathcal O_{\tilde U})=\mathrm H^2(U,\mathcal O_U)$ and this
latter group is zero since we have taken $U$ Stein. Consequently we
have that $\mathrm H^2(\tilde U,\mathcal O_{\tilde U}^*)=0$ and
every $\mathbb P^1$-bundle over $\tilde U$ is the projectivization
of a rank two vector bundle over $\tilde U$.
\end{proof}

The examples below show that the lemma \ref{L:extensao} is no longer
true in dimension greater than two.

\begin{example}\rm
Let $f:\mathbb C^3 \to \mathbb C$ be the function
$f(x,y,z)=x^2+y^2+z^2$ and consider $\mathcal F$ the codimension one
foliation induced by the levels of $f$. If $T\mathcal F$ denotes the
tangent sheaf of $\mathcal F$ then $T\mathcal F$ is a rank $2$
locally free  subsheaf of $T\mathbb C^3$ outside the origin of
$\mathbb C^3$ since at these points $f$ is a local submersion.
Nevertheless, at the origin of $\mathbb C^3$, $T \mathcal F$ is not
locally free, i.e., we cannot write
\[
   df = i_X i_Y dx \wedge dy \wedge dz \, ,
\]
with $X$ and $Y$ germs of holomorphic vector fields at zero. To see
this one has just to observe  that, for arbitrary germs of
holomorphic vector fields $X$ and $Y$, the zero set of $i_X i_Y dx
\wedge dy \wedge dz$ is either empty or has codimension smaller then
two. If $\pi:( \widetilde{\mathbb C^3},D ) \to (\mathbb C^3,0)$
denotes the blow-up of the origin of $\mathbb C^3$ then, as the
reader can check, the tangent sheaf of  $\widetilde {\mathcal F}=
\mathcal \pi^* \mathcal F$ is locally free everywhere. Now the
restriction of $\pi$ to $\widetilde{\mathbb C^3}\setminus \vert D
\vert$ induces an isomorphism of the $\mathbb P^1$-bundles $\mathbb
P(T \widetilde{\mathcal F} \vert_{\widetilde{\mathbb C^3}\setminus
\vert D \vert})$ and $\mathbb P(T {\mathcal F} \vert_{{\mathbb
C^3}\setminus \{0\}})$. Althought the $\mathbb P^1$-bundle $\mathbb
P(T {\mathcal F} \vert_{{\mathbb C^3}\setminus \{0\}})$ does not
extends to a $\mathbb P^1$-bundle over $\mathbb C^3$. \qed
\end{example}

A more geometric version of the previous example has been
communicated to us by C. Ara\'{u}jo. It has appeared several times in
the literature, cf. \cite{Bonavero} and references therein. We
reproduce it here for the reader's convenience.

\begin{example}\rm Let  $p$ be  a point in $\mathbb P^{3}$ and $V$
be the variety of $2$-planes in $\mathbb P^3$ containing $p$.
Consider the variety $X\subset\mathbb P^{3}\times V$ defined as
$$ X=\{(p,\Pi)\in\mathbb P^{3}\times V\,|\,p\in\Pi\}.$$
Consider the natural projection  $\rho: X\to \mathbb P^{3}$. If
$q\neq p$, the fiber over $q$ is the $\mathbb P^{1}$ of planes
containing $p$ and $q$. The fiber  over $p$ is naturally identified
with $V$, thus isomorphic to $\mathbb P^{2}$. If  $\pi:
\widetilde{\mathbb P^3}\to\mathbb P^{3}$ is the blow-up of $p$ then
we have the  following diagram
$$
\xymatrix{ \widetilde{\mathbb P^3} \boxtimes_{\mathbb P^3} X
\ar[d]_{}\ar[r]^{} &
  X\ar[d]^{\rho}
\\ \widetilde{\mathbb P^3} \ar[r]_{\pi}&
{\mathbb P^{3}}}$$ The reader can check that the fibered product
$\widetilde{\mathbb P^3} \boxtimes_{\mathbb P^3} X$ is a $\mathbb
P^1$-bundle over $\widetilde{\mathbb P^3}$ and we are in a situation
analogous to the previous example.\qed
\end{example}

To finish the proof of Theorem \ref{T:normalform} we have to
establish the

\subsection{Unicity of the Normal Form} Let $\mathcal P=(\pi:P\to S,\mathcal H,\sigma:S
\dashrightarrow P)$ and $\mathcal P'=(\pi:P'\to S,\mathcal
H',\sigma:S \dashrightarrow P')$ be two transversely projective
structures in normal for the same foliation $\mathcal F$ and in the
same bimeromorphic equivalence class. Let $\phi:P \dashrightarrow
P'$ be a fibered bimeromorphism. We want to show that $\phi$ is in
fact biholomorphic.

Since both $\mathcal P$ and $\mathcal P'$ are in normal form we have
that  $(\mathcal P)_{\infty}= (\mathcal P')_{\infty}$. Thus for
every $p \in S \setminus \vert (\mathcal P)_{\infty}\vert$ there
exists  a neighboorhood $U$ of $p$ such that  $\mathcal
H\vert_{\pi^{-1}(U)}$ and $\mathcal H'\vert_{\pi'^{-1}(U)}$ are
smooth foliations transverse to the fibers of $\pi$ and $\pi'$,
respectively. If $\phi$ is not holomorphic when restricted to
$\pi^{-1}(U)$ then it most contract some fibers of $\pi$. This would
imply the existence of singular points for $\mathcal
H'\vert_{\pi'^{-1}(U)}$ and consequently contradict our assumptions.
Thus $\phi$ is holomorphic over every $p \in S \setminus \vert
(\mathcal P)_{\infty}\vert$.

Suppose now that $p \in \vert (\mathcal P)_{\infty}\vert$ is a
generic point and that $\Sigma_p$ is germ of curve at $p$ transverse
to $\vert (\mathcal P)_{\infty}\vert$. The restriction of $\phi$ to
$\pi^{-1}(\Sigma)$ (denoted by $\phi_{\Sigma}$) induces a
bimeromorphism of $\mathbb P^1$-bundles over $\Sigma$. Since
$\Sigma$ has dimension one this bimeromorphism can be written as a
composition of elementary transformations. Since $p$ is generic on
the fiber $\pi^{-1}(p)$ we have two of three distinguished points:
one or two singularities of $\mathcal H$ and one point from the
section $\sigma$. But $\phi_{\Sigma}$ must send these points to the
corresponding over the fiber $\pi'^{-1}(p)$. This clearly implies
that $\phi_{\Sigma}$ is holomorphic. From the product structure of
$\mathcal H$ in a neighborhood of $p$, cf. lemma  \ref{L:trivializa}
and remark \ref{R:util} after lemma \ref{L:correcao}, it follows
that $\phi$ is holomorphic in a neighborhood of $\pi^{-1}(p)$.

At this point we have already shown that there exists $Z$, a
codimension two subset of $S$, such that $\phi\vert_{\pi^{-1}(S
\setminus Z)}$ is holomorphic.

Let now $p \in Z$ and $U$ be a neighborhood of $p$ where both $P$
and $P'$ are trivial $\mathbb P^1$-bundles. Thus after restricting
and taking trivializations of both $P$ and $P'$ we have that
$\phi\vert_{\pi^{-1}(U)}$ can be written as
\[
\phi\vert_{\pi^{-1}(U)}(x,[y_1:y_2]) = ( x , [a(x)y_1 + b(x)y_2:
c(x)y_1 + d(x) y_2] ) \, ,
\]
where $a,b,c,d$ are germs of holomorphic functions. But then the
points  $x \in U$ where $\phi$ is not biholomorphic are determined
by the equation $(ad-bd)(x)=0$. Since $(ad-bd)(x)$ is distinct from
zero outside the codimension two set  $Z$ it is distinct from zero
everywhere. Therefore we conclude that $\phi$ is fact biholomorphic
and in this way conclude the prove of the unicity of the normal
form. This also concludes the proof of Theorem \ref{T:normalform}.
\qed

\begin{remark}\rm
To prove the unicity we have not used that $S$ is a surface.
Therefore as long as a normal form exists it is unique no matter the
dimension of the ambient manifold.
\end{remark}

\section{Eccentricity of a Singular Transversely Projective Structure}

\subsection{Foliations on the Projective Plane and  on $\mathbb P^1$-bundles}
The degree of a foliation $\mathcal F$ on $\mathbb P^2$ is defined
as the number of tangencies of $\mathcal F$ with a general line $L$
on $\mathbb P^2$. When $\mathcal F$ has degree $d$ it is defined
through a global holomorphic section of $\mathrm T \mathbb P^2
\otimes \mathcal{O}_{\mathbb P^2}(d-1)$, see \cite[pages
27--28]{Brunella}.

If $\pi:S \to B$ is a $\mathbb P^1$-bundle over a projective curve
$B$, $\mathcal C$ is the foliation tangent to the fibers of $\pi$
and $\mathcal R$ is a Riccati foliation on $S$  then $\mathcal R$ is
defined by a global holomorphic section of $\mathrm TS \otimes
\pi^*(\mathrm T^*B) \otimes \mathcal O_{ S}(\mathrm{tang}(\mathcal
R, \mathcal C))$ \cite[page 57]{Brunella}. If $C\subset S$ is a
reduced curve not $\mathcal R$-invariant then \cite[proposition
2,page 23]{Brunella}
\[
  \deg \left( \pi^*(\mathrm TB)
\otimes \mathcal O_{ S}(-\mathrm{tang}(\mathcal R, \mathcal C)
\right )\vert_C = C^2 - \mathrm{tang}(\mathcal R,C) \, .
\]

With these ingredients at hand we are able to obtain

\subsection{A formula for the Eccentricity: Proof of Proposition \ref{P:ecc}}

Let $L \subset \mathbb P^2$ be a generic line and let $P_L$ be the
restriction of the $\mathbb P^1$-bundle  $\pi:P\to \mathbb P^2$ to
$L$. On $P_L$ we have $\mathcal G$, a Riccati foliation induced by
the restriction of $\mathcal H$, and a curve $C$ corresponding to
$\sigma(L)$. Notice that
\[
 T \mathcal G = (\pi\vert_L)^* \mathcal O_{\mathbb P^1}(2) \otimes
 \mathcal O_{\mathbb P^1}( - (\mathcal P)_{\infty} ) \, .
\]
We also point out that the tangencies between $\mathcal G$ and $C$
are in direct correspondence with the tangencies between $\mathcal
F$ and $L$. Thus
\begin{eqnarray*}
  T \mathcal G \cdot C &=& C\cdot C - \mathrm{tang}(\mathcal G, C)
  \\ &=& - \mathrm{ecc}(\mathcal P) -
  \deg (\mathcal F) \, .
\end{eqnarray*}
Combining this with the expression for $T\mathcal G$ above we obtain
that
\[
  2 - \deg((\mathcal P)_{\infty}) =  - \mathrm{ecc}(\mathcal P) -
  \deg (\mathcal F) \, ,
\]
and the proposition follows. \qed

\subsection{Some Examples}\label{S:exemplos}
Before proceeding let's see some examples of transversely projective
foliations on $\mathbb P^2$ and compute theirs eccentricities using
proposition \ref{P:ecc}.

\begin{example}\rm[Hilbert Modular Foliations on the Projective Plane]
In \cite{LG} some Hilbert Modular Foliations on the Projective Plane
are described. For instance in Theorem 4 of loc. cit. a pair of
foliations $\mathcal H_2$ and $\mathcal H_3$ of degrees $2$ and $4$
is presented. Both foliations admit transversely projective
structures with reduced polar divisor whose support consists of a
rational quintic and a line, cf. \cite{croc2,LG}. For $\mathcal H_2$
the eccentricity is equal to $ 2  = 6 -(2+2)$ while for $\mathcal
H_3$ it is equal to $1 = 6 -(3+2)$. Similarly if one consider  the
pair of  foliations $\mathcal H_5$ and $\mathcal H_9$ presented in
Theorem 2 of loc. cit. then $\mathcal H_5$ has eccentricity $8 = 15
- (5+2)$ and $\mathcal H_9$ has eccentricity $4=15 - (9+2)$. Since
$\mathcal H_5$ is birationally equivalent to $\mathcal H_9$ and
$\mathcal H_2$ is birationally equivalent to $\mathcal H_3$ these
examples show that the eccentricity is not a birational invariant of
transversely projective foliations.
\end{example}

\begin{example}\rm[Riccati Foliations on $\mathbb P^2$]
Let $p\in \mathbb P^2$ be a point and let $\mathcal F$ be a degree
$d$ foliation for which the singular point $p$ has $l(p)=d$. We
recall that $l(p)$ is defined as follows: if $\pi:
(\widetilde{\mathbb P^2},E) \to (\mathbb P^2,p)$ is the blow-up of
$p$ and $\omega$ is a local $1$-form with codimension two singular
set defining $\mathcal F$ then $l(p)$ is the vanishing order of
$\pi^*\omega$ along the exceptional divisor $E$.

When $l(p)=d$ it follows from \cite[page 28 example 3]{Brunella}
that $T \pi^* \mathcal F \cdot \overline L = 0$, where $\overline L$
is the strict transform of a line passing through $p$. From the
discussion in  \cite[page 50--51]{Brunella} it follows that $\pi^*
\mathcal F$ is a Riccati foliation.

For a generic degree $d$ foliation $\mathcal F$ satisfying $l(p)=d$
we will have $d+1$ invariant lines passing through $p$ and no other
invariant algebraic curves. Since $\mathcal F$ is Riccati it will
have a transversely projective structure with exceptional divisor
supported on the $d+1$ $\mathcal F$-invariant lines. For generic
$\mathcal F$ the exceptional divisor will be reduced and with
support equal to the union of  these lines. In this case we will
have that the eccentricity is minus one.
\end{example}

\begin{example}\rm[Brunella's Very Special Foliation]
The very special foliation admits a birational model on $\mathbb
P^2$ where it is induced by the homogeneous $1$-form (cf.
\cite{kod0})
\begin{equation}\label{E:modelo}
   \omega=(- y^2 z - x z^2 + 2xyz) dx + (3 xyz- 3x^2 z) dy + ( x^2 z  -2 x y^2 + x^2 y )dz
   \, .
\end{equation}
It has three invariant curves. The   lines $\{x=0\}$ and $\{z=0\}$
and the rational cubic $\{ x^2 z+x z^2 -3xyz + y^3 =0 \}$. Notice
that the rational cubic has a node at $[1:1:1]$. Moreover
\[
d \omega = \left( \frac{dx}{x} + \frac{dz}{z} + \frac{2}{3} \frac{d(
x^2 z+x z^2 -3xyz + y^3) }{x^2 z+x z^2 -3xyz + y^3 }\right) \wedge
\omega \, .
\]
It can be verified that $\mathcal F$ has a projective structure in
normal form with the polar divisor reduced and with support equal to
the three  $\mathcal F$-invariant curves. Thus the eccentricity of
this projective structure is one.
\end{example}

\subsection{Proof of Proposition \ref{P:dynamics}}

Let $\mathcal F$ be a quasi-minimal singular transversely projective
foliation of $\mathbb P^2$ with transverse structure $\mathcal
P=(\pi:P\to M, \mathcal H, \sigma:M\dashrightarrow P)$ in normal
form. If the monodromy of $\mathcal H$ is non-solvable and not
minimal then there exists a non-algebraic proper closed set
$\mathscr M$ of $P$ formed by a union of leaves and singularities of
$\mathcal H$.

If $L \subset \mathbb P^2$ is  a generic line then
$\mathrm{ecc}(\mathcal P) = -C^2$ where $C = \sigma(L)$. If
$\mathrm{ecc}(\mathcal P)\le0$, i.e., $C^2 \ge 0$  then every  leaf
of $\mathcal G$, the restriction of $\mathcal H$ to $\pi^{-1}(L)$
must intersects $\mathscr M \cap \pi^{-1}(L)$. In the case $C^2>0$
this follows from \cite[Corollary 8.2]{Laurent}. When $C^2=0$ we
have that  $\pi^{-1}(L)=\mathbb P^1 \times \mathbb P^1$ and every
non algebraic leave must intersect every fiber of the {\it
horizontal} fibration(otherwise the restriction of the second
projection to it would be constant).

Therefore for $L$ generic enough $\sigma^* \mathscr M$ is a
non-algebraic proper closed subset of $\mathbb P^2$ invariant under
$\mathcal F$. Thus $\mathcal F$ is not quasi-minimal. This
contradiction implies the result. \qed

\section{The Monodromy Representation}

\subsection{A Local Obstruction}\label{S:obstrucao}

Let $H=\{x_1\cdot x_2 =0 \}$ be the union of the coordinate axis in
$\mathbb C^2$ and $\rho: \pi_1(\mathbb C^2 \setminus H) \to
\mathrm{PSL}(2,\mathbb C)$ a representation.

\begin{prop}
If  $\rho$ is the monodromy representation of a  transversely
projective structure $\mathcal P$ defined in $U$, a neighborhood of
$0 \in \mathbb C^2$, then $\rho$ lifts to $\mathrm{SL}(2,\mathbb
C)$.
\end{prop}
\begin{proof}
We can suppose without loss of generality that $U$ is a polydisc and
that $\mathcal P$ is normal form. Over $U$ every $\mathbb
P^1$-bundle is trivial therefore $\mathcal H$ induces  an integrable
differential $\mathfrak{sl}(2,\mathbb C)$-system on the trivial rank $2$ vector
bundle over $U$, cf. \S \ref{S:singular}. Clearly $\rho$ {\it lifts}
to the monodromy of the $\mathrm{sl}(2,\mathbb C)$-system and the
proposition follows.
\end{proof}

\noindent {\bf A word of warning:} it is not true that the monodromy
of a transversely projective structure $\mathcal P$ always lift to
$\mathrm{SL}(2,\mathbb C)$. For instance we have smooth Riccati
equations over elliptic curves with  monodromy group conjugated to
the abelian group
 $$G=< (z_1:z_2) \mapsto (z_2:z_1); (z_1:z_2)
\mapsto (-z_1:z_2) >.$$

\subsection{Prescribing the monodromy: Proof of Theorem \ref{T:RH}}\label{S:RH}

First we will assume  that $H$ is an hypersurface with smooth
irreducible components and with at most normal crossings
singularities. Instead of working with the projective surface $S$ we
will work with a projective manifold $M$ of  arbitrary dimension
$n$.

\medskip

\noindent{\bf Construction of the $\mathbb P^1$-bundle and of the
foliation.} If $\rho: \pi_1(M \setminus H ) \to \text{SL}(2;\mathbb
C)$ is a representation then it follows from Deligne's work on
Riemann-Hilbert problem \cite{Katz} that there exists $E$, a rank 2
vector bundle over $M$, and a meromorphic flat connection
$$\nabla:E \to E \otimes \Omega^1_M(\log H)$$  with monodromy representation given by $\rho$. From the $\mathbb
C$-linearity  of $\nabla$ we see that its solutions induce $\mathcal
H$, a codimension one  foliation of  $\mathbb P(E)$. If
$\pi_{\mathbb P(E)}: \mathbb P(E) \to M$ denotes the natural
projection then over $\pi_{\mathbb P(E)}^{-1}(M\setminus H)$ the
restriction of $\mathcal H$ is nothing more than suspension of
$[\rho]:\pi_1(M\setminus H) \to \text{PSL}(2,\mathbb C)$ as defined
in \cite[Example 2.8]{croc2}.

Let  $U$ be a sufficiently small open set of $M$ and choose a
trivialization of $E_{\vert U}= U \times \mathbb C^2$ with
coordinates $(x,z_1,z_2) \in U \times \mathbb C \times \mathbb C$.
 Then  for every section $\sigma=(\sigma_1,\sigma_2)$ of $E_{|U}$ we
have that
\[
   \nabla_{|U}(\sigma) = \left( \begin{matrix} d\sigma_1 \\ d\sigma_2 \end{matrix}
   \right) + A \cdot \left( \begin{matrix} \sigma_1 \\ \sigma_2 \end{matrix}
   \right) \,
\]
where \[
 A=\left(\begin{matrix} {\alpha }& \beta\\
\gamma & \delta
\end{matrix}\right)
\]
is two by two matrix with  $\alpha,  \beta, \gamma, \delta \in
\Omega^1_M(\log H)$ satisfying the integrability condition $dA + A
\wedge A=0$. Thus $\nabla = 0$ induces the system
\begin{eqnarray*}
  dz_1 &=&  z_1 \alpha + z_2 \beta  \\
  dz_2 &=&  z_1 \gamma + z_2 \delta \, .
\end{eqnarray*}
Thus the solution of the above differential system are contained in
the leaves of the  foliation defined over $\pi_{\mathbb
P(E)}^{-1}(U)$  by
\[
   \Omega_U = z_1dz_2 - z_2dz_1 - z_2 ^2 \beta  + z_1z_2 (\gamma -\alpha )  + z_1^2 \delta \, .
\]
Clearly the foliations defined in this way patch together to give
$\mathcal H$, a codimension one foliation on $\mathbb P(E)$
transverse to fibers of $\pi$ which are not over $H$.

\medskip

\noindent{\bf Construction of the meromorphic section.} The next
step in the proof of Theorem \ref{T:RH} is to assure the existence
of a {\it generic} meromorphic section of $\mathbb P(E)$. This is
done in the following

\begin{lemma}\label{L:section}
There exists a meromorphic  section
\[
   \sigma:M \dashrightarrow  \mathbb P(E) \, ,
\]
with the following properties:
\begin{enumerate}
\item[(i)] $\sigma$ is generically transversal to  $\mathcal G$;
\item[(ii)] $\overline{\mathrm{sing}(\sigma^*\Omega) \setminus
(\sigma^*\mathrm{sing}(\mathcal G) \cup \mathrm{Ind}(\sigma))}$ has
dimension zero.
\end{enumerate}
\end{lemma}
\begin{proof}
Let $\mathcal L$ be an ample line bundle over $M$. By Serre's
Vanishing Theorem we have that for $k\gg 0$ the following properties
holds:
\begin{enumerate}
\item[(a)] $E\otimes \mathcal L^k$  is generated by global sections;
\item[(b)] for every $x \in M$, $E\otimes \mathcal L^k \otimes m_x$ and $E\otimes \mathcal L^k \otimes
m_x^2$ are also generated by global sections.
\end{enumerate}
Using a variant of the arguments presented in \cite[proposition
 5.1]{Sad} it is possible to settle   that there exists a Zariski  open  $V \subset
\mathrm H^0(M,E\otimes \mathcal L^{\otimes k})$  such that for every
$s \in V$ the zeros locus of $s$ is non-degenerated, of codimension
two, with no irreducible component contained in the support of $H$
and whose image does not contains any irreducible component of
$\mathrm{sing}(\mathcal F)$. We leave the details to the reader.

Let now $\mathcal U=\{U_i\}_{i\in I}$ be a finite covering of $M$ by
Zariski open subsets such that the restrictions of  $E$ and of the
cotangent bundle of $M$ to each $U_i$ are both trivial bundles. For
each $i \in I$ consider
\begin{eqnarray*}
 \Psi_i: U_i \setminus (U_i \cap H) \times \mathrm H^0(M,E\otimes \mathcal L^{\otimes k}) &\to&
 \mathbb C^n \, \\
 (x,s) &\mapsto& s^*\Omega_i(x)
\end{eqnarray*}
where $\Omega_i$ is the $1$-form over $\pi_{\mathbb P(E)}^{-1}(U_i)$
defining $\mathcal G\vert_{U_i}$ and $\Omega^1_{U_i}$ is implicitly
identified with the trivial rank $n$ vector bundle over $U_i$. It
follows from (a) and (b) that for every $x \in M$ there exists
sections in $\mathrm H^0(M,E\otimes \mathcal L^{\otimes k})$ with
prescribed linear part at $p$. Thus if $Z_i = \Psi_i^{-1}(0)$ then
\[
  \dim Z_i = \mathrm h^0(M,E\otimes \mathcal L^{\otimes k} ) .
\]
If  $\rho_i:Z_i \to \mathrm{H}^0(M,E\otimes \mathcal L^{\otimes k})$
is the  natural projection then there exists a Zariski open set $W_i
\subset  \mathrm{H}^0(M,E\otimes \mathcal L^{\otimes k})$  such that
\[
\dim Z_i \le \dim \rho^{-1}(s) + \mathrm h^0(M,E\otimes \mathcal
L^{\otimes k}) .
\]
Thus $\dim \rho^{-1}(s) = 0$ for every $s \in W_i$.

A section $s \in  \left( \bigcap_{i \in I} W_i \right)  \cap V$ will
induce a meromorphic section $\sigma$ of $\mathbb P(E)$ with the
required properties.
\end{proof}

\medskip

\noindent{\bf Unicity.} It remains to prove the unicity in the case
that $\rho$ is non-solvable. We will need the following

\begin{lemma}\label{L:riccati}
Suppose that $\pi:\mathbb P(E) \to M$ has a meromorphic section
$\sigma$  such that the foliation $\mathcal F= \sigma^* \mathcal H$
have non unique transversely projective structure. Then the
monodromy representation of $\mathcal H$ is meta-abelian or there
exists an algebraic curve $C$, a  rational map $\phi:\mathbb P(E)
\dashrightarrow C \times \mathbb P^1$ and Riccati foliation on $C
\times \mathbb P^1$ such that $\mathcal H = \phi^* \mathcal R$.
\end{lemma}
\begin{proof}
After applying a fibered birational map we can assume that $\mathbb
P(E) = M \times \mathbb P^1$ and that $\sigma$ is the
$[1:0]$-section, i.e., if
 \[
  \Omega = z_1dz_2 - z_2 dz_1 + \alpha z_1^2 + \beta z_1\cdot
  z_2 + \gamma z_2^2 \, ,
\]
is the one form defining $\mathcal H$ then $\mathcal F$ is induced
by $\alpha$.

Since $\mathcal F$ has at least two non bimeromorphically
equivalents projective structures then it follows from
\cite[proposition 2.1]{scardua} (see also \cite[lemma 2.20]{croc2})
that there exists a rational function $\ell$ on $M$ such that
\[
   d \alpha = - \frac{d\ell}{2\ell}\wedge \alpha \, .
\]
Thus, after a suitable change of coordinates we can assume that
$\beta = \frac{d\ell}{\ell}.$ From the relation $d\beta = 2 \alpha
\wedge \gamma$ we deduce the existence of a rational function $f \in
k(M)$ such that $\gamma= f \alpha$. Therefore $d\gamma = \beta
\wedge \gamma$ implies that
\[
  \left( \frac{df}{f} - \frac{dl}{l} \right) \wedge \alpha = 0  .
\]
If  $\mathcal F$ does not admit a rational first integral then  $f =
\ell$. Consequently, on the new coordinate system,
\[
  \Omega = z_1dz_2 - z_2 dz_1 + \alpha z_1^2 + \frac{d\ell}{2\ell} z_1\cdot
  z_2 + \ell \alpha  z_2^2 \, .
\]
If  $\Phi(x,[z_1:z_2]) = (x,[z_1:\sqrt{\ell}z_2])$ then we get
\[
  \frac{\Phi^*\Omega}{\sqrt{\ell}} = z_1dz_2 - z_2 dz_1 + (z_1^2 +z_2^2)
  \frac{\alpha}{\sqrt{\ell}} \, \implies d\left( \frac{\Phi^*\Omega}{\sqrt{\ell}(z_1^2 + z_2^2)}\right) = 0 ,
\]
meaning that after a ramified covering the foliation $\mathcal H$ is
induced by a closed $1$-form. Thus $\mathcal H$ has meta-abelian
monodromy.

When $\mathcal F$ admits a rational first integral then it follows
from \cite[Theorem 4.1.(i)]{scardua}(see also \cite[proposition
2.19]{croc2}) that there exists an algebraic curve $C$, a  rational
map $\phi:\mathbb P(E) \dashrightarrow C \times \mathbb P^1$ and
Riccati foliation on $C \times \mathbb P^1$ such that $\mathcal H =
\phi^* \mathcal R$.
\end{proof}

Back to the proof of Theorem \ref{T:RH} we apply lemma
\ref{L:section} to produce a  section $\sigma: M \dashrightarrow
\mathbb P(E)$ generically transversal to $\mathcal H$. If the
transversely projective structure of $\mathcal F=\sigma^* \mathcal
H$ is non unique then  lemma \ref{L:riccati} implies that there
exists an algebraic curve $C$, a  rational map $\phi:\mathbb P(E)
\dashrightarrow C \times \mathbb P^1$ and Riccati foliation on $C
\times \mathbb P^1$ such that $\mathcal H = \phi^* \mathcal R$.
Recall that we are assuming here that $\rho$ is non-solvable.

As we saw in the proof of lemma \ref{L:section} we have a lot of
freedom when choosing $\sigma$. In particular we can suppose that
$\phi \circ \sigma: M \dashrightarrow C \times \mathbb P^1$ is a
dominant rational map. Thus $\mathcal F$ is the pull-back of Riccati
foliation with non-solvable monodromy by a dominant rational map.
The unicity of the transversely projectice structure of $\mathcal F$
follows from \cite[proposition 2.1]{scardua}.

This is sufficient to conclude the proof of Theorem \ref{T:RH} under
the additional assumption on $H$: normal crossing with smooth
ireducible componentes. Notice that up to this point everything
works for projective manifolds of arbitrary dimension.

To conclude  we have just to consider the case where $H$ is an
arbitrary curve on a projective surface $S$. We can proceed as in
the proof of Theorem \ref{T:normalform}, i.e., if we denote by
$p:(\tilde S,\tilde H=p^*H) \to (S,H)$ the desingularization of $H$
then there exists $\tilde \rho: \pi_1(\tilde S, \tilde H) \to
\mathrm{SL}(2,\mathbb C)$ such that $\rho = p_* \tilde \rho$. Thus
we apply the previous arguments over $\tilde S$  and go back to $S$
using lemma \ref{L:extensao}. \qed

\bibliographystyle{amsplain}

\end{document}